\documentclass[12pt]{amsart}
\usepackage{latexsym, amssymb, amsmath, amscd} 
 \usepackage{eucal}

   \newtheorem{rema}{Remark}[section]
    \newtheorem{propo}[rema]{Proposition}

   \newtheorem{theo}[rema]{Theorem}
   \newtheorem{def-theo}[rema]{Definition-Theorem}
 \newtheorem{conj}[rema]{Conjecture}
   \newtheorem{defi}[rema]{Definition}
    \newtheorem{lemma}[rema]{Lemma}
    \newtheorem{corol}[rema]{Corollary}
     \newtheorem{exam}[rema]{Example}
  \newtheorem{rmk}[rema]{Remark}

 \newcommand{\pf}{\noindent{\it Proof:}\hspace{2mm}} 
\newcommand{\epfv}{\hspace{1em}$\Box$\vspace{1em}}

\newcommand{\F}{\mathbb F}
\newcommand{\Q}{\mathbb Q}
\newcommand{\N}{\mathbb N}
\newcommand{\C}{\mathbb C}
\newcommand{\R}{\mathbb R}
\newcommand{\p}{\partial}

\newcommand{\al}{\alpha}
 
\newcommand{\rad}{{\frak r}\hspace{0.3mm}}
\newcommand{\im}{{\rm Im\hspace{0.3mm}}}
\newcommand{\cL}{\mathcal L}

\begin{document}

\title[Mathieu Subspaces of Univariate Polynomial Algebras]
{Mathieu Subspaces of Univariate Polynomial Algebras}
\author{Arno van den Essen and Wenhua Zhao}
\address{A. van den Essen, Department of Mathematics, Radboud University,  
 Nijmegen, Postbus 9010, 6500 GL Nijmegen, The Netherlands. 
Email: essen@math.ru.nl}

\address{W. Zhao, Department of Mathematics, Illinois State University, Normal, 
IL 61790-4520, USA. Email: wzhao@ilstu.edu}

 %\date{\today}

\begin{abstract} We first give a characterization for 
Mathieu subspaces of univariate polynomial algebras over fields 
in terms of their radicals.  We then deduce that
for some classes of classical univariate orthogonal polynomials the
Image Conjecture is true. We also prove two special cases of  
the one-dimensional Image Conjecture for univariate polynomial algebras $A[t]$ over commutative $\Q$-algebras $A$. 
\end{abstract}

   

\thanks{The second author has been partially supported 
by NSA Grant H98230-10-1-0168}

\bibliographystyle{alpha}
    
\maketitle

%\tableofcontents

\renewcommand{\theequation}{\thesection.\arabic{equation}}
\renewcommand{\therema}{\thesection.\arabic{rema}}
\setcounter{equation}{0}
\setcounter{rema}{0}
\setcounter{section}{0} 

\section{\bf Introduction}\label{S1}

{\it The Jacobian Conjecture} has been the subject of much research
over the last seven decades (see \cite{K}, \cite{BCW} and \cite{E1}). 
Various subcases of this still mysterious conjecture
have been verified. Also, several attempts have been made to generalize the conjecture. However, most of these attempts failed. One of these attempts, which is still fully alive, is {\it the Mathieu Conjecture} posed by Olivier Mathieu \cite{M} in 1995. More recently, based on a symmetric reduction of {\it the Jacobian conjecture} obtained independently by 
M. de Bondt and the first author \cite{BE} 
and G. Meng \cite{Me}, the second author gave 
in \cite{HNP} a new equivalent formulation of {\it the Jacobian Conjecture}. This new formulation is 
called {\it the Vanishing Conjecture} and was in turn generalized later in \cite{GVC} to the so-called 
{\it Generalized Vanishing Conjecture}. 
Very recently, the second author posed 
in \cite{IC} an even stronger conjecture, 
namely, {\it the Image Conjecture}.

In a subsequent paper \cite{GIC},  
both the {\it Image} and the {\it Mathieu} 
conjectures were embedded in a general framework, 
by introducing the notion of Mathieu subspaces 
of rings or algebras. This notion forms a 
generalization of the notion of ideals. 
Both {\it the Image Conjecture} and {\it the Mathieu conjecture} 
can be re-expressed as saying that certain 
subspaces are Mathieu subspaces of suitable algebras. 
For a full account and background of these new conjectures 
the reader is referred to the recent survey paper \cite{E2}.

Furthermore, the two-dimensional {\it Jacobian Conjecture}  
can also be expressed directly in terms of 
Mathieu subspaces: it has been shown in \cite{EWZ1} 
that this conjecture is actually equivalent to 
saying that the image of any derivation of the polynomial 
algebra $\C[x,y]$ in two variables,  
whose divergence is zero, is a Mathieu subspace 
of $\C[x,y]$, if its image contains 
the constant polynomial $1$.

The above connections make it clear that it is very desirable to 
get a better understanding of Mathieu subspaces in some general setting. This paper continues the study of these subspaces, as initiated by the second author in the papers \cite{IC}-\cite{MS}.

Now we give a brief description of the content of this paper.
The results discussed here mainly concern Mathieu subspaces of univariate polynomial algebras over fields or UFDs. First, in section \ref{S2} we recall the definition of a Mathieu subspace and that of the radical of an arbitrary subspace. We also describe some of their basic properties, which will be needed later in the paper.  

In section \ref{S3}, for commutative algebras $A$ over fields $k$, 
we recall some results obtained in \cite{MS} concerning 
the radicals of $k$-subspaces of $A$, whose elements are algebraic 
over $k$. In section \ref{S4}, we use these results to derive a 
characterization for Mathieu subspaces of the univariate 
polynomial algebra $k[t]$ over $k$ (see Theorem \ref{T4.1}). 
Some consequences of this characterization are also discussed  
in this section. 

In section \ref{S5}, we use the characterization derived in section 4 to obtain information about 
{\it the Integral Conjectures} 
(see Conjecture \ref{Conj5.1}) posed by the second 
author in \cite{GIC}. One of these conjectures, 
the {\it Image Conjecture} for the univariate Hermite and Jacobi orthogonal polynomials (see Conjecture \ref{Conj5.2}), is studied in detail in section \ref{S6}. 
Finally, in section \ref{S7} we prove two special cases 
of {\it the Image Conjecture} for the univariate polynomial 
algebra $A[t]$ over a commutative $\Q$-algebra $A$ 
(see Theorems \ref{T7.1} and \ref{T7.7}).

\renewcommand{\theequation}{\thesection.\arabic{equation}}
\renewcommand{\therema}{\thesection.\arabic{rema}}
\setcounter{equation}{0}
\setcounter{rema}{0}

\section{\bf Some Basic Properties of Mathieu Subspaces and Their Radicals}
\label{S2}

We first recall the notions of {\it radicals} 
and {\it Mathieu subspaces} introduced respectively 
in \cite{MS} and \cite{GIC} by the second author. Throughout this section $R$ will denote a unital commutative ring and 
$A$ a commutative algebra over $R$. 

For any subset $V$ of $A$, we define the {\it radical} of $V$, 
denoted by $\rad(V)$, to be the subset of all $a\in A$ such that 
$a^m \in V$ for all large $m$. In other words, 
$a$ belongs to $\rad(V)$ if and only if there exists $N\in \N$   
such that $a^m\in V$ for all $m\geq N$. Note that if $V$ is an 
ideal of $A$, then the radical 
$\rad(V)$ of $V$ is just the radical ideal of $V$. 
In particular, $\rad(V)$ itself in this case 
is a radical ideal of $A$. 

\begin{defi}\label{Def-MS}
A subset $V$ of $A$ is said to be a 
{\it Mathieu subspace} of $A$ 
if $V$ is an $R$-subspace or $R$-submodule of $A$ and 
for any $a\in \rad(V)$, the following holds: 
for any $b\in A$, there exists $N\in \N$ 
$($depending on both $a$ and $b$$)$ such
that $a^mb\in V$ for all $m\geq N$.
\end{defi}

Note that the above definition of Mathieu subspaces is 
slightly different from the one given in \cite{GIC} and \cite{MS}. 
But, as shown in Proposition $2.1$ in \cite{MS} these 
two definitions are actually equivalent to each other.  

Note also that any ideal of $A$ is automatically 
a Mathieu subspace of $A$. Therefore, 
the notion of Mathieu subspaces 
provides a natural generalization of the notion of ideals.

To study radicals of arbitrary $R$-subspaces
$V$ of $A$, we consider $I_V$, the {\it largest} ideal  
of $A$ contained in $V$. Since $V$ contains the zero-ideal, 
it is easy to see that $I_V$ always exists and 
is actually the sum of all ideals of 
$A$ contained in $V$. In particular, when $V$ itself 
is an ideal of $A$, we have that $I_V=V$. 

We denote by $\bar A$ the quotient algebra $A/I_V$ 
and by $\pi: A \to \bar A$ the quotient map from $A$ 
to $\bar{A}$. 
We also let $\bar{a}$ and $\overline{V}$ denote $\pi(a)$ and $\pi(V)$, 
respectively. 

The first result of this section provides a different point of view for 
the fact that $I_V$ is the largest ideal contained in $V$.

\begin{lemma}\label{L2.1} 
Let $V$ be an $R$-subspace of $A$. Then for the $R$-subspace $\overline{V}$ of the quotient algebra $\bar A$, 
we have  $I_{\overline{V}}=0$. 
\end{lemma} 
 
\pf Assume otherwise. Then there exists some  
nonzero element $\bar{a}\in \bar{A}$ such that
$\bar a\bar A \subset I_{\overline{V}} \subset \overline{V}$.
Hence $aA+I_V\subset V+I_V=V$, for $I_V\subset V$. Since $\bar{a}$ is nonzero, we have 
$a\not \in I_V$, whence the ideal $aA+I_V$ is strictly larger
than $I_V$. But this contradicts the maximal choice of $I_V$. 
\epfv

The next result shows that radicals correspond nicely under the
algebra homomorphism $\pi: A \to \bar A$.

\begin{propo}\label{P2.2}  
$i)$ $\rad(V)=\pi^{-1}(\,\rad(\pi(V))\,)$. 

$ii)$ $\pi(\rad(V))=\rad(\pi(V))$. 

$iii)$ $\rad(V)$ is an ideal of $A$ 
if and only if $\rad(\pi(V))$ is an ideal of $\bar A$.
\end{propo}

\pf Let $a\in \rad(V)$. Then $a^m\in V$ for all large $m$.
Applying the $R$-algebra homomorphism $\pi$ to $a^m$ gives the 
inclusions $\subset$ in both {\it i}) and {\it ii}).

Now let $\pi(a)\in \rad(\pi(V))$. Then $\pi(a^m)=\pi(a)^m\in \pi(V)$
for all large $m$, say, $\pi(a^m)=\pi(v_m)$ for some $v_m$ in $V$. 
Then $a^m-v_m$ belongs to the kernel of $\pi$, which is equal to $I_V$. 
Since $I_V$ is contained in $V$ and $V$ is an $R$-subspace, 
we have that $a^m\in V$ for all large $m$. Hence we have 
$a\in \rad(V)$, which 
gives respectively the inclusions 
$\supset$ in {\it i}) and {\it ii}). 

Finally, {\it iii}) follows readily from {\it i}) and {\it ii}) 
and the surjectivity of the $R$-algebra homomorphism $\pi$. 
\epfv

Furthermore, Mathieu subspaces also behave nicely under the algebra homomorphism $\pi$, as can be seen from the next proposition 
whose proof is straightforward and is left to the reader 
(or see Proposition $2.7$ in \cite{MS} for a more general statement).

\begin{propo}\label{P2.3}
 $V$ is a Mathieu subspace of $A$ if and
only if $\overline{V}=\pi(V)$ is a Mathieu subspace of $\bar{A}$.
\end{propo}

\begin{rmk}
In the proofs of Propositions \ref{P2.2} 
and \ref{P2.3} given above 
we only use the fact that $I_V$ is an ideal contained in $V$. 
Therefore, the same proof gives the following more general result: 
if $f$ is a surjective $R$-algebra homomorphism from $A$ to $B$, 
and $V$ is an $R$-subspace of $A$ such that the kernel of $f$ 
is contained in $V$, then with $\bar{A}$ replaced by $B$ 
and $\pi$ by $f$, the statements of Propositions \ref{P2.2} 
and \ref{P2.3} also hold.  
\end{rmk} 

From now on, except in section \ref{S7}, the last 
section of the paper, we assume: $k$ is a field, $A$ 
a $k$-algebra and $V$ a $k$-subspace of $A$. 
The importance of Propositions \ref{P2.2} 
and \ref{P2.3} comes from the fact 
that in various situations the quotient algebra $\bar{A}$ turns out
to be algebraic over $k$, i.e., each element of $\bar{A}$ is a root 
of a nonzero univariate polynomial with coefficients in $k$. 
For example, this is the case when $A$ is the univariate polynomial 
algebra $k[t]$ and the ideal $I_V$ is nonzero. We will 
return to this situation in section \ref{S4}. 
It is therefore natural to consider first the case that $A$ 
is algebraic over $k$. This will be done in the next 
section.

\renewcommand{\theequation}{\thesection.\arabic{equation}}
\renewcommand{\therema}{\thesection.\arabic{rema}}
\setcounter{equation}{0}
\setcounter{rema}{0}

\section{\bf Mathieu Subspaces with Algebraic Radicals}
\label{S3}

Suppose that $V$ is a Mathieu subspace of $A$. Then what can be said about
the structure of the radical $\rad(V)$ of $V$?

In this section we discuss this question under the additional
assumption that all elements of the radical $\rad(V)$ are algebraic
over $k$. The next result, first obtained in \cite{MS},  
asserts that the Mathieu subspaces, whose radicals are 
algebraic over $k$, are completely characterized by their radicals.
More precisely, by Theorems 4.10 and 4.12 in \cite{MS}  
we have the following theorem.

\begin{theo}\label{T3.1}
If all elements of $\rad(V)$ are algebraic over $k$, 
then $V$ is a Mathieu subspace of $A$
if and only if $\rad(V)$ is an ideal in $A$. 
In this case, we also have $\rad(V)=\rad(I_V)$.
\end{theo}

\begin{corol} \label{C3.2}
Let $V$, $\bar{A}=A/I_V$ and $\overline{V}$ as before. 
Assume that all elements of $\rad(\overline{V})$ are algebraic over $k$. 
Then $V$ is a Mathieu subspace of $A$ if and only of 
$\rad(V)$ is an ideal of $A$. In this case, we also have 
$\rad(V)=\rad(I_V)$. 
\end{corol}

\pf First, the equivalence follows readily from 
Proposition \ref{P2.3}, Theorem \ref{T3.1} and Proposition \ref{P2.2} 
{\it iii}). Second, by Theorem \ref{T3.1} and Lemma \ref{L2.1} we have 
$\rad(\overline{V})=\rad(I_{\overline{V}})=
\rad(\overline{0})$. 
It then follows from Proposition \ref{P2.2} {\it i}) 
that $\rad(V)=\pi^{-1}(\rad(\overline{0}))$.
Since the latter set is also equal to $\rad(I_V)$ 
(for $ {\rm Ker\,\,} \pi=I_V$), 
we obtain that $\rad(V)=\rad(I_V)$. 
\epfv 

\renewcommand{\theequation}{\thesection.\arabic{equation}}
\renewcommand{\therema}{\thesection.\arabic{rema}}
\setcounter{equation}{0}
\setcounter{rema}{0} 

\section{\bf A Characterization of Mathieu Subspaces of 
Univariate Polynomial Algebras over Fields} \label{S4}

Throughout this section $k$ is a field and $A=k[t]$, the univariate polynomial
algebra over $k$. Although, apart from the constant polynomials, 
the algebra $A$ has no elements which are algebraic over $k$, 
we will show that in this
case Mathieu subspaces can also 
be characterized by their radicals. 

\begin{theo}\label{T4.1}
For any $k$-subspace $V$ of $A$, 
$V$ is a Mathieu subspace of $A$ 
if and only if $\rad(V)=\rad(I_V)$.
\end{theo}

\pf First, if the ideal $I_V$ is nonzero, 
$\bar{A}=A/I_V$ is finite dimensional and hence algebraic over $k$. 
Then the theorem follows immediately from Corollary \ref{C3.2}.
So from now on we assume $I_V=0$. 

If $\rad(V)=\rad(I_V)$, then $\rad(V)=\rad(0)=0$, since $A$ has no zero-divisors. By Definition \ref{Def-MS} $V$ obviously 
is a Mathieu subspace of $A$.

Conversely, assume that $V$ is a Mathieu subspace of $A$ (with $I_V=0$). 
We must show that $\rad(V)=\{0\}$. 
Assume the contrary and let $a\in A=k[t]$ be a nonzero element 
in $\rad(V)$. Then $a^m\in V$ for all large $m$. 
If the polynomial $a$ has degree $0$, it is a nonzero constant, 
and hence so are all $a^m$ ($m\ge 1$). 
It follows that $V$ contains $1$. 
Since $V$ is a Mathieu subspace of $A$, it is easy to check 
(or see Lemma $4.5$ in \cite{GIC}) that in this case we have $V=A$. 
Hence we have $I_V=A$, which is a contradiction 
since $I_V=0$. Therefore, we have 
$d\!:=\deg a\geq 1$. 

Since $V$ is a Mathieu subspace of $A=k[t]$, 
there exists $N\ge 1$ such that $t^i a^m\in V$ for all 
$0\leq i\leq d-1$ and all $m\geq N$. In particular, we have 
\begin{align}  
a^m h\in V  \mbox{ \it for all } m\geq N \mbox{ \it and all } h\in k[t] 
\mbox{ \it with } \deg  h\leq d-1.  \label{T4.1-pe1} 
\end{align} 

Now, let $m\geq N$. Since $I_V=0$, the nonzero ideal $a^mA$ cannot be contained 
in $V$, so there exists a nonzero $b_m\in A$ of the lowest degree 
such that $a^m b_m \notin V$. Let $b_{m_0}$ have the smallest degree amongst all the $b_m$ $(m\ge N)$. Note that by the property in 
Eq.\,(\ref{T4.1-pe1}) we have $\deg b_{m_0}\geq d$. 
Furthermore, by the Euclidean division 
there exist $q$ and $r$ in $A$ such that 
$\deg r\leq d-1$ and 
\begin{align}
b_{m_0}=qa+r. \label{T4.1-pe2}  
\end{align}

Since $\deg b_{m_0}\geq d> \deg r$, it follows that 
$\deg qa=\deg b_{m_0}$. Since $\deg a\geq 1$, 
we deduce that $\deg q<\deg b_{m_0}$. 
Furthermore,  by multiplying $a^{m_0}$ to 
Eq.\,(\ref{T4.1-pe2}), we have  
\begin{align}\label{T4.1-pe3}
b_{m_0}a^{m_0}=qa^{m_0+1}+ra^{m_0}.
\end{align}

Now, by the choices of $b_m$'s with $m=m_0$,  
the left hand side of the equation above does not belong to $V$. 
By the property in Eq.\,(\ref{T4.1-pe1}), 
we see that $ra^{m_0}\in V$. Hence, $qa^{m_0+1}\notin V$. 
But then, by the choices of $b_m$'s with $m=m_0+1$, 
we obtain that $\deg b_{m_0+1} \leq \deg q < \deg b_{m_0}$,   
which contradicts the choice of $b_{m_0}$, 
since $b_{m_0}$ has the least degree among all 
the $b_m$ $(m\ge N)$. 
\epfv

One immediate consequence of Theorem \ref{T4.1} is 
the following necessary condition for a $k$-subspace $V$ 
of the univariate polynomial algebra $A=k[t]$ to be 
a Mathieu subspace of $k[t]$.   

\begin{corol}\label{C4.2}
Let $V$ be a $k$-subspace of $A=k[t]$. Assume that $V$ 
is a Mathieu subspace of $A$. Then $\rad (V)$ 
is a radical ideal of $A$.  
\end{corol}

Another consequence of Theorem \ref{T4.1} is  
that for Mathieu subspaces $V$ of $A=k[t]$, the integer  
$N$ in Definition \ref{Def-MS} can actually be chosen 
in a way that depends only on the element $a\in \rad(a)$, 
i.e., $N$ can be chosen to be independent with  
the element $b\in A$ in Definition \ref{Def-MS}. 
More precisely, we have the following corollary.      

\begin{corol}
Let $V$ be a $k$-subspace of $A=k[t]$. 
Then $V$ is a Mathieu subspace of $A$ if and only if for any 
$a\in \rad(V)$, there exists $N\ge 1$ such that 
for all $b\in A$, we have $a^m b\in V$ 
for all $m\ge N$. 
\end{corol}

\pf The $(\Leftarrow)$ part is obvious. To show the 
$(\Rightarrow)$, assume that $V$ is a Mathieu subspace of $A$ 
and let $a\in \rad(V)$. Then by Theorem \ref{T4.1}, 
we have $a\in \rad(I_V)$. Since $I_V$ is an ideal of 
$A$, there exists $N\ge 1$ such that 
$a^m\in I_V$ and hence $a^mA\subset I_V$ for all $m\ge N$. 
Consequently, for any $b\in A$, we have $a^mb \in I_V\subset V$ 
for all $m\ge N$, whence the corollary follows. 
\epfv
 
Finally, we conclude this section with the following example, 
which shows that Theorem \ref{T4.1} does not hold in general 
for polynomial algebras in two or more variables.  

\begin{exam} 
Let $B=\C[x,y]$, $D=x\p_x-y\p_y$ and $V=\im D=D(\C[x, y])$, 
i.e., the image of the derivation $D$. 
Then as shown in the proof of Lemma $3.4$ 
in \cite{EWZ1}, $V$ is a Mathieu subspace of $B$ 
and $\rad(V)=W_1 \cup W_2$, where $W_1$ is the $\C$-span of all 
the monomials $x^iy^j$ with $i<j$ and $W_2$ is the $\C$-span 
of all the monomials $x^iy^j$ with $i>j$. So $\rad(V)$ is
not even a $\C$-subspace of $B$ and hence, not an ideal 
of $B$. 
\end{exam}

\renewcommand{\theequation}{\thesection.\arabic{equation}}
\renewcommand{\therema}{\thesection.\arabic{rema}}
\setcounter{equation}{0}
\setcounter{rema}{0}

\section{\bf Integral and Image Conjectures in Dimension One}
\label{S5}

The aim of this section is to show how the results of the previous section
can be used to obtain some new results concerning several conjectures
of the second author posed in \cite{GIC}. To keep this paper 
as much self-contained as possible, here we briefly recall these 
conjectures for the one-dimensional case.

\begin{conj}\label{Conj5.1} $(${\bf Integral Conjecture}$)$ 
Let $B\subset\R$ be an open subset and
$\sigma$ a positive measure such that $\int_B g(t)d\sigma$ exists and is finite 
for all $g\in\C[t]$. Set 
\begin{align}\label{Def-VB} 
V_B(\sigma)\!:=\Big\{f\in\C[t]\,\Big|\,\int_B f \, d\sigma=0\Big\}. 
\end{align}
Then $V_B(\sigma)$ is a Mathieu
subspace of $\C[t]$. 
\end{conj}

In \cite{GIC} this conjecture is proved for several special cases. One of them is the case when $\sigma$ is an {\em atomic measure} supported at finitely many points $r_i$ in $B$, i.e., 
$\sigma(r_i)>0$ for each $i$ and, for any subset $U\subset B$, 
$\sigma(U)$ is the sum of $\sigma(r_i)$ over all the  
$r_i$'s that are contained in $U$.  
It is proved in Proposition $3.11$, \cite{GIC} that 
$V_B(\sigma)$ in this case is a Mathieu subspace 
of $\C[t]$ by showing that $\rad(V_B(\sigma))$ is 
the ideal of all the polynomials 
vanishing at all $r_i$'s. 

\medskip

To describe the second conjecture, we need to recall some 
results on univariate orthogonal polynomials. 
Let $B$ be an open interval of $\R$ and 
$w(t)$ a so-called {\em weight function} on $B$,
i.e.,  $w(t)$ is non-negative over $B$ and its integral over $B$  
is finite and positive. To such a pair $(B,w)$ one can associate 
a {\it Hermitian inner product} on $\C[t]$ by defining
\begin{align}\label{Def-Biform} 
\langle f,g\rangle=\int_B f(t)\bar g(t)\, w(t)dt,
\end{align}
where $\bar g(t)$ is the complex conjugate of the polynomial of $g(t)$, i.e., the polynomial obtained by taking the complex conjugates of 
the coefficients of $g(t)$. 

Applying the Gram-Schmidt process to the standard basis 
$1,t,t^2,\cdots$ of $\C[t]$,   
we obtain a set of orthogonal polynomials of $\C[t]$. 
Making the following special choices for $B$ and $w(t)$, 
we get the following {\em classical univariate orthogonal polynomials}. 
\begin{enumerate}
  \item[$1)$]  
The {\em Hermite polynomials}: $B=\R$ and $w(t)=e^{-t^2}$. 
  
\item[$2)$] The {\em generalized Laguerre polynomials}: $B=(0,\infty)$ and $w(t)=
t^{\alpha}e^{-t}$ with $\alpha>-1$. 

\item[$3)$] The {\em Jacobi polynomials}: $B=(-1,1)$ and
$w(t)=(1-t)^{\alpha}(1+t)^{\beta}$ with $\alpha,\beta>-1$.
\end{enumerate}

In each of the three cases above, 
we can define the differential operator
\begin{align}\label{Def-DLambda}
\Lambda=w^{-1}\circ\partial_t\circ w.
\end{align}

\noindent This gives respectively the following related operators:
\begin{align}\label{Three-DLambda}
\partial_t-2t,\,\,\,\partial_t+(\alpha t^{-1}-1),\,\,\,\partial_t-\alpha(1-t)^{-1}+\beta(1+t)^{-1}.
\end{align}

\medskip 

In \cite{GIC} the second author makes the following conjecture.

\begin{conj}\label{Conj5.2}
$(${\bf Image Conjecture for classical orthogonal polynomials}$)$ 
Let $\Lambda$ be as defined in Eq.\,$($\ref{Three-DLambda}$)$ with $\alpha,\,\beta>-1$.   
Set 
\begin{align}\label{Conj5.2-e1}
\im'\Lambda:=\C[t]\cap\Lambda(\C[t]).
\end{align}
Then $\im'\Lambda$ is a Mathieu subspace of $\C[t]$.
\end{conj}

Furthermore, the following result has also been proved 
in Lemma $2.5$ $c)$ and Proposition $3.3$ in \cite{GIC}.

\begin{propo}\label{P5.3}
With the same notations as in Conjectures \ref{Conj5.1} and \ref{Conj5.2},  
we have  
\begin{enumerate}
  \item[$i)$] $1\in \im'\Lambda$ if and only if $\im'\Lambda=\C[t]$; 
  \item[$ii)$] if $1\not \in \im'\Lambda$, then 
   $\im'\Lambda=V_B(\sigma)$.  
Consequently, in this case Conjecture \ref{Conj5.1} 
      holds for the pair $(B, \sigma)$ with $d\sigma=wdt$ 
        if and only if Conjecture \ref{Conj5.2} holds for 
             the related differential operator $\Lambda$ 
in Eq.\,$($\ref{Three-DLambda}$)$. 
\end{enumerate}
\end{propo}

Next, we show that Conjectures \ref{Conj5.1} and \ref{Conj5.2} 
are actually respectively equivalent to the following 
two formally stronger conjectures. 

\begin{conj} \label{Conj5.3} 
$(${\bf Strong Integral Conjecture}$)$ 
With the same notations as in Conjecture \ref{Conj5.1}, 
assume that $\sigma$ is not an atomic measure supported 
at finitely many points.  
Then $\rad(V_B(\sigma))=\{0\}$.
\end{conj} 

In other words, the conjecture above claims that 
when the measure $\sigma$ is not an atomic 
measure supported at finitely many points,  
the only polynomial $f$ with $\int_B f^m d\sigma=0$ 
$(m\ge 1)$ should be the zero polynomial.  

\begin{conj} \label{Conj5.4} 
$(${\bf Strong Image Conjecture for classical orthogonal polynomials}$)$ 
With the same notations as in 
Conjecture \ref{Conj5.2}, assume 
$1\notin \im'\Lambda$. Then 
$\rad(\im'\Lambda)=\{0\}$. 
\end{conj} 

\begin{theo}\label{T5.3} 
$i)$ The Integral Conjecture $($Conjecture \ref{Conj5.1}$)$ 
is equivalent to the Strong Integral Conjecture 
$($Conjecture \ref{Conj5.3}$)$. 

$ii)$ The Image Conjecture for classical univariate orthogonal
polynomials $($Conjecture \ref{Conj5.2}$)$  
is equivalent to the Strong Image Conjecture 
for classical univariate orthogonal polynomials 
$($Conjecture \ref{Conj5.4}$)$.  
\end{theo}

\pf {\it i}) Note first that if $\sigma$ is an atomic measure supported at 
finitely many points, then by Proposition $3.3$ in \cite{GIC}, 
Conjecture \ref{Conj5.1} holds. When $\sigma$ is not an atomic measure 
supported at finitely many points, it is easy to see 
from Definition \ref{Def-MS} 
that Conjecture \ref{Conj5.1} follows directly  
from Conjecture \ref{Conj5.3}. Therefore, 
in any case the $(\Leftarrow)$ part of statement 
{\it i}) holds.

To show the $(\Rightarrow)$ part, put $V=V_B(\sigma)$. 
We claim that the largest ideal $I_V$ of $\C[t]$ 
contained in $V$ is equal to $0$, 
which combining with Theorem \ref{T4.1} will imply 
Conjecture \ref{Conj5.3}, for the polynomial algebra 
$\C[t]$ has no zero-divisors. 

Assume the contrary and let $0\ne f\in I_V$.  
Then $g:=\bar{f}f \in I_V \subset V$ (where $\bar{f}$ denotes 
the complex conjugate of $f$), whence $g\in V$. 
Then by definition of $V$, the integral of $g$ over $B$ 
is equal to zero. 
On the other hand, since $g$ is continues and positive over $B$ 
(except at the finitely many zeroes of $f$ in $B$) and $\sigma$ is 
not an atomic measure supported at finitely many points, 
the integral of $g$ over $B$ is positive, which is a contradiction. 

{\it ii}) Note first that if $1\in \im'\Lambda$, then 
by Proposition \ref{P5.3} {\it i}) we have 
$\im'\Lambda=\C[t]$, which is obviously a Mathieu subspace of 
$\C[t]$. If $1\not\in \im'\Lambda$, then 
it is easy to see from Definition \ref{Def-MS} and 
Proposition \ref{P5.3} {\it ii}) that Conjecture \ref{Conj5.2} follows directly from Conjecture \ref{Conj5.4}. Therefore, 
in any case the $(\Leftarrow)$ part of statement 
{\it ii}) holds.

To show the $(\Rightarrow)$ part, assume Conjecture \ref{Conj5.2} 
and $1\not \in \im'\Lambda$. Then by Proposition \ref{P5.3} {\it ii}), 
we have $\im'\Lambda=V_B(\sigma)$, where $d\sigma=wdt$. So the radicals of the two subspaces are equal. Since by our hypothesis 
$V_B(\sigma)$ is a Mathieu subspace of $\C[t]$, 
by Theorem \ref{T4.1} we have 
$\rad(\im'\Lambda)=\rad(V_B(\sigma))=\rad(I_V)$, where 
$V=V_B(\sigma)$ as above. But, as shown in the proof 
of statement {\it i}) above, we also have $I_V=0$ and 
$\rad(I_V)=\{0\}$. Hence $\rad(\im'\Lambda)=\{0\}$, 
as desired.  
\epfv

\renewcommand{\theequation}{\thesection.\arabic{equation}}
\renewcommand{\therema}{\thesection.\arabic{rema}}
\setcounter{equation}{0}
\setcounter{rema}{0}

\section{\bf Some Cases of The Strong Image Conjecture 
for the Hermite and Generalized Laguerre Polynomials}
\label{S6}

In this section, we prove some cases of the Strong Image 
Conjecture ({\bf SIC}), Conjecture \ref{Conj5.4}, and also 
of the following conjecture, which is the one dimensional 
case of {\it the Image Conjecture} posed by 
the second author in \cite{IC}. 

\begin{conj}\label{1d-DIC} 
Let $A$ be any $\Q$-algebra, $c\in A$ and $a(t)\in A[t]$. 
Set $D\!:=c \p_t-a(t)$ and $\im D\!:=D(A[t])$. Then $\im D$  
is a Mathieu subspace of $A[t]$. 
\end{conj} 

The main result of this section is the following theorem.  

\begin{theo}\label{T6.2}
Let $d\in \N$ and $\al\in \Q$ such that 
$\al \not \in -(1+(d+1)\N)$ and 
$(d,\al)\neq (0,0)$. 
Let $D=\p_t+\al t^{-1}-t^d$. Then $\rad(\im'D)=\{0\}$, 
where $\im'D\!:=\C[t]\cap D(\C[t])$. 
\end{theo}

One consequence of the theorem above is the following 
corollary on the case of the {\bf SIC} for the Hermite and generalized Laguerre polynomials. 

\begin{corol}\label{T6.2-C1}
The {\bf SIC} holds for the Hermite polynomials and 
the generalized Laguerre polynomials with 
$\alpha\in\Q$ $($and $\alpha>-1$$)$. 
\end{corol}

\pf For the Hermite polynomial case, we make the variable change 
$t=\sqrt{2}s$. Since $\p_s+\al s^{-1}-2s=\sqrt{2}(\p_t+\al t^{-1}-t)$, by Theorem \ref{T6.2} with $d=1$ we see that the {\bf SIC} holds in this case.  

The generalized Laguerre polynomial case follows from Theorem \ref{T6.2} by taking $d=0$, 
in case $\al\neq 0$. If $\al=0$, 
then $D=\p_t-1$, which is an invertible map with the inverse 
$D^{-1}=-\sum_{k\ge 0}\p_t^{k}$. Hence $\im'D=\C[t]$.
In particular, $1\in \im' D$ and the condition 
of the {\bf SIC} does not apply in this case. 
But, since $\C[t]$ is obviously a Mathieu subspace 
of $\C[t]$, we see that Conjecture \ref{Conj5.2}  
still holds in this case. 
\epfv

Another consequence of Theorem \ref{T6.2} is the 
following special case of Conjecture \ref{1d-DIC}. 

\begin{corol}\label{T6.2-C2}
Conjecture \ref{1d-DIC} holds 
for the case that $A=\C$ and 
$a(t)= \lambda t^d$ for all $\lambda\in \C$ 
and $d\ge 0$.
\end{corol}

\pf Note that in this case $D=c\p_t-\lambda t^d$ is a differential operator of $\C[t]$ itself. Hence we have $\im' D=\im D$. 
If $c=0$, then $D=-\lambda t^d$ and $\im D$ is the ideal of $\C[t]$ generated by $\lambda t^d$. If $c\ne 0$ but $\lambda =0$, 
then $D=c\p_t$, whence $\im D=\C[t]$. Therefore, in either 
of these two cases, $\im D$ is an ideal of $\C[t]$ and hence, 
is also a Mathieu subspace of $\C[t]$. 

Assume that both $c$ and $\lambda$ are nonzero.  
If $d=0$, then $D=c\p_t-\lambda$, which is invertible with 
the inverse $D^{-1}=-\sum_{i\ge 0}\lambda^{-i+1}c^i\p^i_t$. 
Hence $\im D=\C[t]$, which is obviously 
a Mathieu subspace of $\C[t]$. Therefore, Conjecture \ref{1d-DIC} 
also holds in this case.  

So, we may further $d\ge 1$. Let $\beta\in \C$ such 
that $\beta^{d+1}=c/\lambda$. 
By changing of variables $t=\beta s$, it is easy to check 
that $D=c\beta^{-1}(\p_s-s^d)$.   
Then by the fact that $\im' D=\im D$ (as mentioned above), 
Conjecture \ref{1d-DIC} follows immediately 
from Theorem \ref{T6.2} (for the different operator $\p_s-s^d$) 
and Definition \ref{Def-MS}.  
\epfv

Next, we give a proof for Theorem \ref{T6.2} starting with the following observations. 

\medskip

First, since $D t^n=(n+\al)t^{n-1}-t^{d+n}$ for all $n\geq 1$, 
we have  
\begin{align} 
t^{n+d}  \equiv (n+\al)\, t^{n-1} \quad (\mbox{\rm mod } \im'D). \label{T6.2-pe1}  
\end{align} 
Applying the relation above repeatedly, 
it follows that for any $k\ge d+1$, we have 
\begin{align} 
t^k  \equiv c_k \, t^i \quad (\mbox{\rm mod } \im'D), \label{T6.2-pe1b}  
\end{align} 
where $0\le i\le d$ with  
$i \equiv k \mod (d+1)$ and $c_k$ is given by 
\begin{align*}
c_k=\big(k-(d+1)+1+\alpha \big)
\big(k-2(d+1)+1+\alpha\big)\cdots \big(i+1+\alpha\big). 
\end{align*}  
\noindent Therefore, we can define a 
$\C$-linear map $\cL$ from $\C[t]$ to the vector subspace $S$ of polynomials of degree $\leq d$, by setting $\cL(t^k)=c_k t^i$ for all $k\ge d+1$, and $\cL(t^j)=t^j$ for all $0\le j\le d$. 
Then $\cL$ has the property 
that for any $h\in \C[t]$, $h\in \im'D$ if and only if
$\cL(h)\in S \cap \im'D$. 

Furthermore, we also define 
the $\C$-linear functional $\cL_0:\C[t]\to \C$ by 
setting $\cL_0(h)\!:=\cL(h)(0)$ for all $h\in \C[t]$,  i.e., we set $\cL_0(h)$ to be the constant term of 
the polynomial $\cL(h)\in S$. 

\begin{lemma}
Let $d\in \N$, $\alpha \in \C$ and $D=\p_t+\al t^{-1}-t^d$. 
Assume $(d, \al)\ne (0, 0)$.  
Then we have
 \begin{align}
\im'D \subseteq {\rm Ker\,\,} \cL_0. \label{T6.2-pe2} 
\end{align}
\end{lemma}
\pf First, if $\al\neq 0$, then    
$\al t^{-1}-t^d=D\cdot 1 \not \in \im' D$. It is easy to see that 
in this case every nonzero element of $\im'D$ has degree 
at least $d+1$, whence $S\cap \im'D=\{0\}$. 
If $\al=0$, then $d\ge 1$ by the assumption 
of the theorem. Since in this case $t^d=D(-1)\in \im' D$, 
it is easy to see that $S \cap \im'D=\C \, t^d$, 
the one-dimensional subspace spanned by $t^d$.   

Hence, in any case we have that $f(0)=0$ for all 
$f\in S\cap \im' D$. Since for any $h\in \C[t]$, 
$h \in \im'D$ if and only if $\cL(h)\in S\cap \im'D$, 
we see that the lemma follows.
\epfv
  
Second, the $\C$-linear functional $\cL_0$ can be described more explicitly by 
the next lemma. But, we need first to fix the 
following notation: for any 
$\alpha\in \C$ and positive integers $q, n\in \N$, we set 
 \begin{align}\label{qnn!} 
[qn,n]_{\al}!:=((q-1)n+1+\al)((q-2)n+1+\al)\cdots(1+\al). 
\end{align} 
Furthermore, for convenience we also set (for the case $q=0$)   
\begin{align}\label{0nn!} 
[0,\, n]_{\al}!:=1 
\end{align} 
for all $\alpha\in \C$ and integers $n\ge 1$. 

Note that if $\al\not \in -(1+(d+1)\N)$, then it follows from Eq.\,(\ref{qnn!}) 
that $[q(d+1),\, (d+1)]_{\al}\ne 0$ for all $q\in \N$.

\begin{lemma}\label{L6.3}
Let $d$ and $\alpha$ be as in Theorem \ref{T6.2}. 
Then for any $0\leq i\leq d$ and $q\in \N$, 
we have 
\begin{align}\label{L6.3-e1}
\cL_0\big(t^{q(d+1)+i}\big)=
\begin{cases} 0, &\mbox{ if } i>0,\\
[q(d+1),\, d+1]_{\al}! &\mbox{ if } i=0.
\end{cases}
\end{align} 
\end{lemma} 
\pf First, by definition of $\cL$ we have that
$\cL(t^{q(d+1)+i})=c_it^i$ for some $c_i\in\C$. 
Consequently, if $i>0$, the constant term of $\cL(t^{q(d+1)+i})$ 
is zero, which gives the first case of Eq.\,(\ref{L6.3-e1}). 

Second, by Eq.\,(\ref{T6.2-pe1}) with $n=(q-1)(d+1)+1$ 
we have 
\begin{align*}
t^{q(d+1)} \equiv \big((q-1)(d+1)+1+\al\big)\, t^{(q-1)(d+1)} \,\,\, (\mbox{\rm mod } \im'D).
\end{align*} 
Then by applying the induction on $q$, the second case 
of Eq.\,(\ref{L6.3-e1}) also follows. 
\epfv

\noindent\underline{\it Proof of Theorem \ref{T6.2}\,}: 
Assume otherwise. We  
fix a nonzero $f(t)\in \rad(\im'D)$ and derive 
a contradiction as follows. 

First, we write 
\begin{align}   
f=c_st^s+c_{s+1}t^{s+1}+\cdots+c_N t^N \label{T6.2-pe3}  
\end{align}
for some integers $s \le N$ and $c_i\in \C$ $(s\le i \le N)$ 
with $c_s,\, c_N\ne 0$. 

Note that we may obviously assume that $c_s=1$. Moreover, 
we may also assume that $s\geq 1$. To see this, suppose that $s=0$ and 
that we have already proved the $s\geq 1$ case. Let $m_0\ge 1$ be such that 
$f^m\in \im'D$ for all $m\geq m_0$ and set $g:=f^{m_0}-f^{m_0+1}$. Then it is easy to see 
that $g(0)=0$ (for $f(0)=1$) and $g^m\in \im'D$ for all $m\geq 1$. 
By our assumption on the $s\geq 1$ case, we have $g=0$. 
Since $f\ne 0$, we deduce that $f=1$, whence $1\in \im'D$.
But this is obviously impossible, since $(d,\al)\neq (0,0)$.

Next, by a similar reduction used by M. Boyarchenko in his unpublished proof (but see \cite{FPYZ}) for the case of 
Conjecture \ref{Conj5.3} with $B=[0, 1]\subset \R$ 
and $d\sigma=dt$ (see also \cite{EWZ2} or \cite{E2} 
for a similar reduction), we may also assume 
that all the coefficients $c_i$'s of $f$ in Eq.\,(\ref{T6.2-pe3}) 
belong to some algebraic number field $K$. 
Then by a well-known result in algebraic number theory (e.g., see \cite{Weiss}), 
we know that for each prime $p\ge 2$, there exists at least one extension 
of the $p$-valuation $v_p(\cdot)$ of $\Q$ to $K$, which we will still 
denote by $v_p(\cdot)$, such that for all but finitely 
many prime numbers $p$, we have $v_p(c_i)\geq 0$ for all 
$s\le i\le N$. 

Now we consider $f^{m(d+1)}$ $(m\ge 1)$. Note that for all large $m$ 
this element belongs to $\im'D$, and hence by Eq.\,(\ref{T6.2-pe2}) it also belongs to ${\rm Ker\,\,} \cL_0$, i.e., $\cL_0(f^{m(d+1)})=0$ for all large $m$. From our reductions on $f$, we have $s\ge 1$ and $c_s=1$. Hence  
we also have  
\begin{align}\label{T6.2-pe4}
f^{m(d+1)}=t^{sm(d+1)}+\sum_{k\geq sm(d+1)+1} \phi_k t^k,
\end{align}
where the $\phi_k$'s are polynomials in the $c_i$'s 
with integer coefficients. 

Now apply $\cL_0$ to Eq.\,(\ref{T6.2-pe4}) and observe that 
by Lemma \ref{L6.3} the only powers of $t$ 
on the right hand side of the equation, which contribute to 
$\cL_0(f^{m(d+1)})$, are the powers $t^k$ with $k$ divisible by $d+1$. So for all $m\gg 0$ we have 
\begin{align}\label{T6.2-pe5}
\cL_0\big(t^{sm(d+1)}\big)+\sum_{i\geq 1}\phi_{(sm+i)(d+1)} 
\cL_0\big(t^{(sm+i)(d+1)}\big)=0.
\end{align} 
Then by Eq.\,(\ref{L6.3-e1}), we get 
\begin{align*}
[sm(d+1),\, d+1]_{\al}!+\sum_{i\geq 1}[(sm+i)(d+1),\, d+1]_{\al}!\,
\phi_{(sm+i)(d+1)}=0. 
\end{align*}

\noindent Dividing by $[sm(d+1),\, d+1]_{\al}!$ from the equation above, we get 
\begin{align*}
1+\sum_{i\geq 1}b_i\, \phi_{(sm+i)(d+1)}=0,   
\end{align*}
\begin{align}
\sum_{i\geq 1}b_i\, \phi_{(sm+i)(d+1)} = -1, 
\label{T6.2-pe6}  
\end{align}
where the coefficients $b_i$'s are given by 
\begin{align}\label{T6.2-pe6b}
 b_i=\big( (sm+i-1)(d+1)+1+\al \big)\cdots \big(sm(d+1)+1+\al\big). 
\end{align}

Assume first that $\al\neq 0$. Note also that $\al\neq -1$, 
since by assumption $\al \not \in -(1+(d+1)\N)$. 
Write $\al=r/q$ for some nonzero integers such that 
$q\ge 1$ and $\gcd(r,q)=1$. Then
\begin{align}\label{T6.2-pe7}
 b_i=\big( (sm+i-1)q(d+1)+q+r \big)\cdots \big( smq(d+1)+q+r \big)\big/q^i. 
\end{align}

Observe that $q$ and the numerator of $b_i$   
have no common factor. We claim that 
for any large enough $m$, there exists  
a prime number $p_m$ which divides 
the numerators of all the $b_i$'s.  
 
Observe first that $smq(d+1)+q+r$ divides the numerators of all the $b_i$'s and that $q+r \neq 0$ (since $\al\neq -1$). 
Let $s_0=\gcd(s(d+1),q+r)$, $s(d+1)=s_0 s_*$ and $q+r=s_0 h$.    
Hence $\gcd(s_*, h)=1$. 
Then we have 
\begin{align}\label{T6.2-pe8}
smq(d+1)+q+r=(s_0s_*q)m+s_0h=s_0\big((s_*q)m+h\big).
\end{align}

Since $\gcd(s_*q,h)=1$, it follows from Dirichlet's prime number theorem (e.g., 
see Theorem $66$ and Corollary $4.1$, 
p.\,$297$ in \cite{FT}) that 
there exists infinitely many $m\ge 1$ such that 
$p_m\!:=(s_*q)m+h$ is a prime 
number. Note that by Eqs.\,(\ref{T6.2-pe7}) and 
(\ref{T6.2-pe8}) any such a prime number 
$p_m$ divides the numerators 
of all the $b_i$'s in Eq.\,(\ref{T6.2-pe6}). 

Now choose and fix any such large enough $m\ge 1$, and write $p_m$ as $p$ for short,  
such that an extension $v_p(\cdot)$ of 
$p$-valuation of $\Q$ to the number field $K$ 
satisfies that $v_p(c_i)\geq 0$ for all 
$s\le i\le N$. Then we also have $v_p(\phi_j)\geq 0$ for all $\phi_j$'s 
in Eq.\,(\ref{T6.2-pe6}). 
Since $v_p(b_i)>0$ for all $b_i$'s in 
Eq.\,(\ref{T6.2-pe6}), 
it follows that the $v_p$-valuation of 
the left hand side of Eq.\,(\ref{T6.2-pe6}) 
is positive. Consequently, from 
Eq.\,(\ref{T6.2-pe6}) we have 
$v_p(-1)>0$. But, this is 
a contradiction since $v_p(-1)=0$. 

Finally, we consider the case $\al=0$. 
Note that in this case by Eq.\,(\ref{T6.2-pe6b}),  
all the nonzero $b_i$'s in  
Eq.\,(\ref{T6.2-pe6}) are positive integers 
that are divisible by $sm(d+1)+1$. 
Then by Dirichlet's prime number theorem again, 
there exist infinitely many $m\ge 1$ 
such that $p_m\!:=sm(d+1)+1$ is a prime. 
Applying the same argument as above 
we will get a contradiction again. 
Therefore, the theorem follows.  
\epfv

Next we consider the case when the condition 
$\al\notin -(1+(d+1)\N)$ in 
Theorem \ref{T6.2} fails. 

\begin{theo}\label{T6.65} 
Let $d\in \N$, $\al\in -(1+(d+1)\N)$ and $D=\p_t+\al t^{-1}-t^d$. 
Then we have 
\begin{enumerate} 
\item[${\it i})$] The statement of Theorem \ref{T6.2} 
does not hold. More precisely, we have $t^{d+1}\in \rad(\im' D)$. 
\item[${\it ii})$] If $d\geq 1$, then 
$(t^{d+1})^m t=t^{(d+1)m+1} \notin\im' D$ for all
$m\geq 0$. So $\im' D$ is not a Mathieu subspace of 
$\C[t]$. 
\item[${\it iii})$] If $d=0$, then $\im' D$ is a Mathieu subspace of $\C[t]$ and $\rad(\im' D)=t\,\C[t]$.
\end{enumerate}  
\end{theo}

\pf {\it i}) Let $\al=-(1+q(d+1))$ for some $q\in\N$. Then for each $m\geq 0$, by Eq.\,(\ref{T6.2-pe1}) with $n=(m+q)(d+1)+1$ we get 
\begin{align} 
t^{(m+q+1)(d+1)}\equiv m(d+1)t^{(m+q)(d+1)} \quad (\mbox{\rm mod } \im'D). \label{T6.65-pe1}
\end{align} 
%$$D(t^{(m+q)(d+1)+1})=m(d+1)t^{(m+q)(d+1)}-t^{(q+m+1)(d+1)}.$$
%\noindent So $t^{(q+m+1)(d+1)}\equiv m(d+1)t^{(q+m)(d+1)} (\mbox{ mod } \im' D)$ for all $m\geq 0$. 
In particular, by choosing $m=0$ we see that 
$t^{(q+1)(d+1)}\in \im' D$. Then from this fact and 
Eq.\,(\ref{T6.65-pe1}), it is easy to see that for any $k\geq q+1$, 
we have $t^{k(d+1)}\in \im' D$, whence $t^{d+1}\in \rad(\im' D)$.

{\it ii}) In a similar way, for each $m\geq 1$, by Eq.\,(\ref{T6.2-pe1}) with $n=(m-1)(d+1)+2$ we get 
\begin{align} 
 t^{m(d+1)+1}\equiv \big((m-1-q)(d+1)+1\big)t^{(m-1)(d+1)+1}   (\mbox{\rm mod } \im'D). \label{T6.65-pe2}
\end{align} 
  
Since $d\geq 1$ and the factor appearing on the right hand side 
of the equation above is equivalent $1$ modulo $d+1$, we see that 
this factor can not be equal to zero. 
It then follows by applying Eq.\,(\ref{T6.65-pe2}) repeatedly that 
$t^{(d+1)m+1}\equiv c\,t$ $(\mbox{\rm mod } \im'D)$ for some nonzero 
$c\in\C$. On the other hand, under the assumptions $d\ge 1$ and  
$\al\in -(1+(d+1)\N)$ it is easy to verify directly 
that $\al\neq 0$ and 
$t \notin \im' D$. Hence statement 
{\it ii}) follows. 

{\it iii}) Since $d=0$, we have $\alpha=-(1+q)$. Then 
for any $m\ge 0$, by Eq.\,(\ref{T6.2-pe1}) with $d=0$ and $n=m+1$ 
we have 
\begin{align} 
 t^{m+1}\equiv (m-q) t^m \quad   (\mbox{\rm mod } \im'D). \label{T6.65-pe3}
\end{align} 
Applying this equation repeatedly we see that  
$t^n\in\im' D$ for all $n\geq q+1$. Hence,  
$t^{q+1}\C[t]\subseteq \im' D$ 
and $t\,\C[t]\subseteq \rad(\im' D)$.  

On the other hand, since $\al\le -1$ in this case, 
it is easy to check that $1 \notin \im' D$. 
Then from Corollary $7.12$ in \cite{MS} 
with $\mathcal A=\C[t]$ and the maximal ideal 
$\frak m=(t)$ it follows that    
$\im' D$ is indeed a Mathieu subspace 
of $\C[t]$.  Moreover, by Corollary \ref{C4.2} 
we see that $\rad(\im' D)$ is an ideal of $\C[t]$. 
Since $t\,\C[t]\subseteq \rad(\im' D)$ 
(as pointed out above) and $1\notin \rad(\im' D)$ 
(for $1\notin \im' D$), we also have  
$\rad(\im' D)=t\,\C[t]$. Hence, statement 
{\it iii}) follows. 
\epfv

To conclude this section we point out that the {\bf SIC} also holds for the following special Jacobi polynomials. But, let's first recall  the following results related with the univariate Jacobi polynomials. 

\begin{theo}\label{T6.7}
Let $B=(-1, 1)\subset\R$, $w(t)=(1-t)^\alpha(1+t)^\beta$ $(\alpha, \beta>-1)$, 
$d\sigma=w(t)dt$ and $V_B(\sigma)$ as 
in Eq.\,$($\ref{Def-VB}$)$. Then we have  
\begin{enumerate}
  \item[{\it i}$)$] If $\al, \beta \in \N$, then $\rad(V_B(\sigma))=0$.
\item[{\it ii}$)$] If $\al=\beta=\lambda-\frac{1}{2}$ and 
$\lambda\in\frac{1}{2}\N$, then $\rad(V_B(\sigma))=0$. 
\end{enumerate} 
\end{theo}

\pf  {\it i}) follows from Theorem 3.4 and Corollary 3.5 in \cite{P}. {\it ii}) is exactly the content of Proposition 4.2 in \cite{FPYZ}.
\epfv

From the theorem above, we see that the {\it Strong Integral Conjecture} (Conjecture \ref{Conj5.3}) holds for all the {\it Gegenbauer polynomials} (i.e., the Jacobi polynomials with 
$\al=\beta=\lambda-\frac{1}{2}$) 
with $\lambda\in\frac{1}{2}\N$. 
Hence, {\it Strong Integral Conjecture} also holds for the {\em Chebyshev polynomials of 
the first and the second kind}, i.e., the 
{\it Gegenbauer polynomials} with 
$\lambda=0,1$, respectively, 
and also for the {\it Legendre polynomials}, 
i.e., the {\it Gegenbauer polynomials} with   
$\lambda=\frac{1}{2}$.

To see whether or not the {\bf SIC} also holds for the differential operator $D=\p_t-\al(1-t)^{-1}+\beta(1+t)^{-1}$ $(\al, \beta>-1)$ related with  
the Jacobi orthogonal polynomials, we need first 
to show the following lemma.    

\begin{lemma}\label{L6.8}
Let $D=\p_t-\al(1-t)^{-1}+\beta(1+t)^{-1}$ with
$\al,\beta >-1$. Then $1\in \im' D$ 
if and only if $\al=0$ or $\beta=0$.
\end{lemma}

\pf The $(\Leftarrow)$ part can be easily checked. 
For example, if $\al\neq 0$ and $\beta=0$, we have 
$D(t-1)=(1+\al)\neq 0$, for $\al>-1$, whence  
$1\in \im' D$. The other case is similar.  
   
Conversely, let $1=D(h)$ for some
$h\in\C[t]$. Assume otherwise, i.e., 
both $\alpha$ and $\beta$ are nonzero. 
Then multiplying the equation $D(h)=1$ by $1-t^2$ we obtain
\begin{align}
(1-t^2)\partial_t h-\alpha (1+t)h+\beta (1-t)h=1-t^2. 
\label{L6.8-pe1} %\,\,\,\,(*)$$
\end{align}

\noindent It follows from the equation above 
that both $1+t$ and $1-t$ divide $h$, so $h=(1-t^2)g$
for some $g\in \C[t]$. Substituting this equality into 
Eq.\,(\ref{L6.8-pe1}) and then dividing by
$1-t^2$ gives  
\begin{align}
(1-t^2)\partial_t g-(2t)g-\alpha (1+t)g
+\beta (1-t)g=1. \label{L6.8-pe2} %\,\,\,\,(**)
\end{align}

Now, write $g=c_d t^d+$ {\it lower order terms}, 
with $0\neq c_d\in\C$. Then by comparing 
the coefficients of $t^{d+1}$ on both sides of 
Eq.\,(\ref{L6.8-pe2}), we get
\begin{align*}
-(d+2)-(\alpha+\beta)=0. %\label{L6.8-pe3}
\end{align*}

\noindent It then follows that $\alpha+\beta=-(d+2)\leq -2$, 
which is a contradiction since
both $\alpha$ and $\beta$ are greater than $-1$.
\epfv

%\pf $(\Leftarrow)$ can be checked directly 
%(or see  Proposition $3.8$ in \cite{GIC}). 
%To show the $(\Rightarrow)$ part, we assume contrary and let $h(t)\in \C[t]$ such that 
%$1=Dh(t)$. Note that for the Jacobi polynomials, we have $B=(-1, 1)\subset \R$ with the {\it weight function} $w(t)=(1-t)^{\alpha}(1+t)^{\beta}$.   
%Since $D$ as a differential operator is the 
%same as $w(t)^{-1}\circ \p_t\circ w(t)$ 
%(as pointed in Eqs.\,(\ref{Def-DLambda}) and 
%(\ref{Three-DLambda}), we have 
%\begin{align}
%1=w(t)^{-1}\p_t\big( h(t) w(t)\big).
%\end{align}

%Now consider the integral 
%\begin{align*}
%\int_{-1}^1 1\cdot w(t)dt &=\int_{-1}^1 (D h) \, w(t)dt =
%\int_{-1}^1 w(t)^{-1}\p_t\big( h(t) w(t)\big) \, w(t)dt \\
%&=\int_{-1}^1  \p_t\Big( h(t) w(t)\big) \, dt =
% h(t) w(t)\big|_{-1}^1=0,
%\end{align*}
%where the last equality follows from the condition that 
%$\alpha, \beta>0$, which implies $w(-1)=w(1)=0$.

%But on the other hand, since $w(t)$ is a continuous 
%function with $w(t)>0$ for all $t\in (-1, 1)$, 
%we also have  
%$\int_{-1}^1 1\cdot w(t)dt=\int_{-1}^1 w(t)dt>0$, which is a contradiction. 
\medskip

Now, from Proposition \ref{P5.3}, Theorem \ref{T6.7} and Lemma \ref{L6.8} we immediately get the 
following corollary. 

\begin{corol}\label{C6.9}
Let $D=\p_t-\al(1-t)^{-1}+\beta(1+t)^{-1}$ such that 
$\al$ and $\beta$ are not both zero. Then 
the following statements hold. 
 \begin{enumerate}
   \item[{\it i}$)$] If $\al, \, \beta \in \N$, then the {\bf SIC} holds for $D$, i.e., $\rad(\im' D)=0$.

\item[{\it ii}$)$] If $\al=\beta=\lambda-\frac{1}{2}$ with $\lambda \in \frac{1}{2}\N$, then $\rad(\im' D)=0$. 
\end{enumerate} 
\end{corol} 
  
In particular, the {\bf SIC} holds for the differential operators related with the {\it Gegenbauer polynomials} 
with $\lambda\in\frac{1}{2}\N$ and 
$\lambda\ge 1$, and also for 
differential operators related with 
the {\em Chebyshev polynomials of 
the first and the second kind}. 
For the {\em Legendre polynomials}, 
i.e., the case $\lambda=\frac12$, 
we have $D=\p_t$ and $1\in \im'D=\C[t]$, 
the condition of the {\bf SIC} is not 
satisfied. But,    
the {\it Strong Integral conjecture} 
obviously still holds in this case 
(as already pointed before).

\renewcommand{\theequation}{\thesection.\arabic{equation}}
\renewcommand{\therema}{\thesection.\arabic{rema}}
\setcounter{equation}{0}
\setcounter{rema}{0}

\section{\bf The One Dimensional Image Conjecture over a Commutative Ring}
\label{S7} 

In this section, we study more special cases of Conjecture \ref{1d-DIC} for commutative $\Q$-algebras $A$. It has been shown in Theorem $2.8$ in \cite{EWZ2} that the conjecture holds under the  assumptions that $a$ is a non-zero-divisor and  
that $aA$ is a radical ideal of $A$.  
We first show in the next theorem that when $A$ 
is a UFD, this radical hypothesis can actually be dropped. 

\begin{theo}\label{T7.1} 
Let $A$ be a UFD and a $\Q$-algebra. Let $a\in A$ and  
$D\!:=\p_t-a$.  Then $\im D$ is a Mathieu subspace of 
$A[t]$. 
\end{theo} 

\begin{rmk}\label{R7.5}
In case that $A$ is an $\F_p$-algebra 
$($not necessarily a UFD$)$, by Theorem $2.2$ in  \cite{EWZ2} the theorem above also holds provided the element $a\in A$ is not a zero-divisor.  Furthermore, it has also been shown  
in Corollary $2.6$ in \cite{EWZ2} that 
for all $f\in A[t]$ with $f^p\in \im D$, 
the $p$-th power of each coefficient of $f$ belongs  to $aA$. From this result and the fact that each monomial $a^pt^n=D^p((-1)^pt^n)$ 
belongs to $\im D$, one deduces easily that $\rad(\im D)=\rad(aA[t])$. In particular, $\rad(\im D)$ in this case is a $($radical$)$ ideal of $A[t]$.    
\end{rmk}

To prove Theorem \ref{T7.1}, note that the case $a=0$ is trivial. 
So, from now on we assume $a\ne 0$ and derive first 
a lemma as follows.    

We denote by $\cL$ the $A$-linear map from $A[t]$ 
to $A$ defined by $\cL(t^n)=n!$ for all $n\geq 0$ and
by $\rad(a)$ the radical ideal of $aA$. 
Since $a$ is not a zero-divisor 
of $A$, it is easy to see that:
\begin{align}
\mbox{\it 
for any $b\in A$,\,\,\,$b\in \im D$ if and only if $b\in aA$ } \tag{$\ast$}
\end{align}

\begin{lemma} \label{L7.2}  With the setting above, the following statements hold.
\begin{enumerate}
\item[{\it i})] for any $n\geq 0$, $a^n t^n \equiv n!\,\, ({\rm mod} \,\,\im D)$. Furthermore, $a^{n+1}t^n A$ is contained in $\im D$ for all $n\geq 0$. 

\item[{\it ii})] Let $f=p(at)$ for some $p(t)\in A[t]$. 
  Then $f\in \im D$ if and only if $\cL(p)\in aA$. 

\item[{\it iii})] Let $f$ be as in {\it ii}$)$. Then $f\in \rad(\im D)$ if and only if all coefficients of $p(t)$ belong to $\rad(a)$. 
\end{enumerate}
\end{lemma}

\pf {\it i}) The first statement follows by induction on 
$n\ge 1$ from the equality $D(a^{n-1}t^n)=na^{n-1}t^{n-1}-a^nt^n$. The second statement follows from the first one and the equivalence in $(\ast)$.

{\it ii}) By {\it i}) and the definition of the $A$-linear map 
$\cL$, we have that $f=p(at)\in \im D$ if and only if 
$\cL(p)\in \im D$. Then the statement follows immediately 
from the equivalence in $(\ast)$.

{\it iii}) Write $p(t)=\sum c_i t^i$ with $c_i\in A$. It follows from {\it ii}) that $f\in \rad(\im D)$ if and only if $\cL(p^m)\in aA$ for all large $m$. Clearly, if all
$c_i$'s belong to $\rad(a)$, then for all large $m$, 
$\cL(p^m)$ belongs to $aA$, whence 
$f\in \rad(\im D)$. 

Conversely, assume that $\cL(p^m)$ belongs to $aA$ for all large $m$ and suppose that some $c_i$'s, say $c_{i_0}$, does not belong to $\rad(a)$. Then by the fact that $\rad(a)$ is the intersection of all prime ideals containing $aA$, there exists a prime ideal 
$\frak{p}$ of $A$ which contains $a$ but not $c_{i_0}$. 

Now, let $R\!:=A/\frak{p}$ and $K$ the field 
of fractions of $R$. Define $\cL_1: K[t]\to K$ 
to be the $K$-linear map such that $\cL_1(t^n)=n!$ 
for all $n\ge 0$. Then by viewing the reduced 
polynomial $\bar p\in R[t]$ inside $K[t]$, 
we have $\cL_1(\bar p^m)=0$ for all large $m$. 
Then by Lefschetz's principle and Theorem 4.9 in \cite{EWZ2}, we have 
$\bar p=0$, i.e., all the coefficients $c_i$'s 
of $p$ lie in the prime ideal $\frak{p}$, which is a 
contradiction. $\Box$ 

\begin{corol}\label{C7.3}
Let $f=p(at)$ for some $p\in A[t]$. If $f\in \rad(\im D)$,
then for every $g\in A[t]$, we have $g f^m \in \im D$ for all large $m$. 
\end{corol}

\pf Let $g\in A[t]$ with degree $d\ge 0$. By Lemma \ref{L7.2} {\it iii}), all coefficients of $p(t)$ belong to 
$\rad(a)$. Hence there exists $N\ge 1$ such that all coefficients of $p(t)^N$ belong to $aA$. 
Consequently, all coefficients of $p(t)^{N(d+1)}$ 
belong to $a^{d+1}A$ and the same holds for $p(t)^m$ 
whenever $m\geq N(d+1)$. Then for each 
$m\geq N(d+1)$ and each $i\geq 0$, 
the coefficient of $t^i$ in $f^m=p(at)^m$ 
belongs to the ideal $a^{i+d+1}A$, 
whence for each $j\ge 0$ the coefficient 
of $t^j$ in $f^mg$ lies in the ideal $a^{j+1}A$.     
It then follows from Lemma \ref{L7.2} {\it i}) 
that $gf^m\in \im D$ for all $m\geq N(d+1)$.
\epfv

\noindent\underline{\it Proof of Theorem \ref{T7.1}\,}: 
First, if $a$ is a unit in $A$, then since 
$\p_t$ is locally nilpotent on $A[t]$, 
the operator $D$ as an $A$-linear map is invertible 
with the inverse map given by 
$D^{-1}=-\sum_{i\ge 0} a^{-i-1}\p_t^i$, 
whence $\im D=A[t]$ and obviously is a Mathieu 
subspace of $A[t]$. 

So we may assume that $aA$ is a proper ideal $A$. 
Hence, so is the ideal $I:=\cap_{i=0}^{\infty} a^iA$. 
Then for each $c\in A\backslash I$, there exists a unique integer $n$ such that $c\in a^nA$ but 
$c\notin a^{n+1} A$. We define $v_a(c)\!:=n$ 
in this case, and set $v_a(c)\!:=\infty$ for 
all $c\in I$. Furthermore, we extend $v_a(\cdot)$ to 
$A[t]$ by setting $v_a(c\,t^i):=v_a(c)-i$ for all 
$c\in A$ and $i\ge 0$.

Now let $f\in \rad(\im D)$ and write $f=\sum_{i=0}^d c_i t^i$ with $c_i\in A$. If all the coefficients $c_i$'s belong to $I$, then 
$f$ certainly can be written in the form $p(at)$ for some 
$p(t)\in A[t]$. Then for any $g\in A[t]$, by Corollary \ref{C7.3} we have $f^mg\in \im D$ when $m\gg 0$. 
So we may assume that not all coefficients of 
$f$ belong to $I$. Let $s(f)$ be the minimum of all 
$v_a(c_it^i)$ $(0\le i\le d)$. If $s(f)\geq 0$, 
then by Corollary \ref{C7.3} again we are done. 

So assume $s(f)\leq -1$. Then $f_1:=a^{-s(f)}f$ also belongs 
to $\rad(\im D)$ and $s(f_1)=0$. 
Write $f_1=\sum b_i t^i$. Then it follows 
that each $b_i$ is of the form $b_i=a^i d_i$ 
for some $d_i\in A$, and furthermore
that $d_i\notin aA$ for some $i$. 
Let $p(t)=\sum d_i t^i$. Then $f_1(t)=p(at)$ and   
by Lemma \ref{L7.2} {\it iii}) $d_i\in \rad(a)$ 
for each $i$. 

By Lemma \ref{L7.4} below, there exist $u\in A$ and
$\tilde{d}_i\in A$ such that $ud_i=\tilde{d}_i a$ 
for all $i$, and $\tilde{d}_i\notin \rad(a)$ 
for some $i$. Observe that since $-s(f)\geq 1$, 
the polynomial $a^{-s(f)-1}f$ 
belongs to $\rad(\im D)$. Hence, so does the polynomial  
$f_2:=u a^{-s(f)-1}f$. 

Now we consider    
\begin{align*}
f_2(t)&=a^{-1}u f_1(t)=
a^{-1} u \sum  b_i t^i =a^{-1} u \sum  d_i a^i t^i\\ 
&=\sum  a^{-1} (u d_i) a^i t^i =\sum \tilde{d}_ia^it^i.  
\end{align*} 
Hence we have $f_2(t)=q(at)$   
with $q(t)=\sum \tilde{d}_i t^i \in A[t]$.
Then applying Lemma \ref{L7.2} to $f_2\in \rad(\im D)$, 
we see that all $\tilde{d}_i\in \rad(a)$. But this is a contradiction, since as pointed out above $\tilde{d}_i \notin \rad(a)$ for some $i$. $\Box$

\begin{lemma} \label{L7.4}
Let $d_1,\cdots,d_n$ all be in $\rad(a)$ and $d_i\notin aA$
for some $i$. Then there exist elements $u\in A$ and $\tilde{d}_i\in A$, such that
$ud_i=\tilde{d}_ia$ for all $i$, and $\tilde{d}_i\notin \rad(a)$ for some $i$.
\end{lemma}

\pf Let $b$ be the greatest common divisor of $a$ and
all $d_i$'s. Since some $d_i\notin aA$, there exists an  irreducible factor $p$ of $a$ such that 
its multiplicity in $b$ is smaller than its multiplicity in $a$. Let $u=a/b$ and $\tilde{d}_i=d_i/b$ for all $i$. Then $ud_i=a\tilde{d}_i$ for each $i$.
Furthermore, if $\tilde{d}_i\in \rad(a)$ for all $i$, 
then $p$ divides each $\tilde{d}_i$
and hence $pb$ is a common divisor of $a$ 
and all $d_i$'s. But this contradicts the definition 
of $b$.
\epfv 

Next, we consider the special case of 
Conjecture \ref{1d-DIC} under the further 
assumption that $1\in \im D$. 
Note that if Conjecture \ref{1d-DIC} 
holds under this assumption, 
then we will have $\im D=A[t]$. 
This is because of the general 
fact that any Mathieu subspace 
$V$ of an algebra $B$ with $1\in V$ 
must be the whole algebra $B$, as first 
noticed in Lemma $4.5$ in \cite{GIC}. 
Indeed, since $1^m=1\in V$ 
for all $m\ge 1$, then for any $u\in B$, 
by taking large enough $m$, we have 
$u=u1^m\in V$, whence $V=B$. 

\begin{theo}\label{T7.7} 
Let $A$ be any $\Q$-algebra, $a(t)\in A[t]$ and $D\!:=c\p_t-a(t)$ for some $c\in A$. Assume $1\in \im D$. Then $\im D=A[t]$.
\end{theo}
 
Before we prove this theorem we make some preparations. First, if
$a(t)=0$, the hypothesis that $1\in \im D$ implies that $c$ is a unit in $A$, which in turn implies that $\im D=A[t]$ (since $A$ is a $\Q$-algebra). So we may assume that $a(t)\neq 0$. Furthermore, 
we can also reduce to the case that $A$ is Noetherian. More precisely,  let $1=D(h)$ for some $h(t)\in A[t]$ and $A_0$ the Noetherian $\Q$-subalgebra of $A$ generated by $c$ and  coefficients of $a(t)$ and 
$h(t)$. Then $D$ as an operator on $A_0[t]$ (by restriction) also has the property $1\in D(A_0[t])=\im D_{A_0[t]}$. So, if we can prove the Noetherian case, it follows that $A_0[t]=D(A_0[t])$. 
In particular, each monomial $t^n$ belongs to 
$D(A_0[t])\subset D(A[t])$, whence $\im D=A[t]$.

To prove the Noetherian case, we will make a reduction to the domain case by using the following result 
in commutative algebra. 

\begin{lemma}\label{L7.8}
For any Noetherian ring $A$, the zero ideal $(0)$ is a 
product of finitely many prime ideals of $A$. 
\end{lemma}

\pf Since $A$ is Noetherian, its nilradical $\frak{n}$
can be written as $\frak{p}_1\cap\cdots\cap\frak{p}_r$ for some 
prime ideals $\frak{p}_i$. 
Since $\frak{n}^e=0$ for some $e\geq 1$, 
it follows that 
\begin{align*}
(0)\subset(\frak{p}_1\cdots\frak{p}_r)^e
\subset(\frak{p}_1\cap\cdots\cap\frak{p}_r)^e
=\frak{n}^e=(0).
\end{align*}
Hence we have $(0)= \frak{p}_1^e\frak{p}_2^e \cdots\frak{p}_r^e$.  
\epfv

\noindent\underline{\it Proof of Theorem \ref{T7.7}\,}:  
As pointed above, we may assume that $A$ is Noetherian and
that $a(t)$ is nonzero. By Lemma \ref{L7.8} we write  
$(0)=\frak{p}_1\cdots\frak{p}_s$ for
some prime ideals $\frak{p}_i$. 

We next show that we may also assume that $A$ is a domain.  
Namely, let $\bar{A}\!:=A/\frak{p}_1$ and $\overline{D}$ 
be the induced operator on $\bar{A}[t]$. 
If we can prove the domain case, it follows
that $\bar{A}[t]=\overline{D}(\bar{A}[t])$. 
Then for any $f\in A[t]$, we have 
$f=D(b)+ \sum_i p_i a_i$ for some $b\in A[t]$,
$p_i\in\frak{p}_1$ and $a_i\in A[t]$. 
By a similar result with $\frak{p}_1$ 
replaced by $\frak{p}_2$, we obtain that
$a_i=D(b_i)+ \sum p_{ij} a_{ij}$ for some $b_i\in A[t]$,  
$p_{ij} \in \frak{p}_2$ and $a_{ij}\in A[t]$. 
Combining these two results we get 
\begin{align*}
f\in D(A[t])+\frak{p}_1\frak{p}_2A[t]. 
\end{align*}
Repeating this argument we finally find that 
\begin{align*}
f\in D(A[t])+\frak{p}_1\cdots\frak{p}_sA[t]=D(A[t]). 
\end{align*}

So it remains to prove the case that $A$ is a domain. However, in this case by comparing the degrees of both sides of the equation 
$1=(c\p_t-a(t))h$, it is easy to see that 
$\deg a(t)=\deg h(t)=0$, i.e., both $a(t)$ and $h(t)$ 
actually belong to $A$. 
It then follows that $a(t)$ is a unit in $A$. 
So we may assume that $a(t)=1$. Since $c\p_t$ is locally nilpotent on $A[t]$, it is easy to see that $D$ is invertible with the inverse $D^{-1}=-\sum_{i\ge 0} c^i\p_t^i$. Hence $\im D=A[t]$, and the theorem follows.
\epfv

\end{document}